\documentclass[11pt]{article}

\usepackage[margin=1in]{geometry}
\usepackage{amsmath,amssymb,amsthm,mathtools}
\usepackage{enumitem}
\usepackage{microtype}
\usepackage{needspace}
\usepackage[colorlinks=true,linkcolor=blue,citecolor=blue,urlcolor=blue]{hyperref}

\hypersetup{
  pdftitle={Spectral Response Geometry for Liouville Brownian Motion},
  pdfauthor={Chunhao Cai},
  pdfsubject={Liouville Brownian motion, spectral response geometry, isospectral deformations},
  pdfkeywords={Liouville Brownian motion, Gaussian multiplicative chaos, spectral response, eigenfunction products, isospectral deformations}
}

\numberwithin{equation}{section}
\newtheorem{theorem}{Theorem}[section]
\newtheorem{proposition}[theorem]{Proposition}
\newtheorem{lemma}[theorem]{Lemma}
\newtheorem{corollary}[theorem]{Corollary}
\theoremstyle{definition}

\theoremstyle{remark}

\newcommand{\R}{\mathbb R}

\newcommand{\Prob}{\mathbb P}
\newcommand{\Exp}{\mathbb E}
\newcommand{\Hc}{\mathcal H}
\newcommand{\E}{\mathcal E}

\newcommand{\abs}[1]{\left|#1\right|}
\newcommand{\inner}[2]{\left\langle #1,#2\right\rangle}
\newcommand{\Tr}{\operatorname{Tr}}
\newcommand{\supp}{\operatorname{supp}}
\newcommand{\spanop}{\operatorname{span}}
\newcommand{\op}{\mathrm{op}}
\newcommand{\GammaE}{\Gamma_{\E}}
\newcommand{\Green}{G_D}
\newcommand{\Lzero}{L_0}

\newcommand{\so}{\mathrm{so}}
\newcommand{\SO}{\mathrm{SO}}
\newcommand{\Sym}{\operatorname{Sym}}

\title{Spectral Response Geometry for Liouville Brownian Motion}
\author{Chunhao Cai\\[0.4em]
\small School of Mathematics (Zhuhai), Sun Yat-sen University, Zhuhai, Guangdong, China\\
\small\href{mailto:caichh9@mail.sysu.edu.cn}{\texttt{caichh9@mail.sysu.edu.cn}}}
\date{}

\begin{document}
\maketitle

\begin{abstract}
For Dirichlet Liouville Brownian motion associated with a subcritical Gaussian
multiplicative chaos measure on a bounded planar domain, we study the response
of finite spectral blocks under Cameron--Martin perturbations of the underlying
Gaussian free field. For each fixed $0<\gamma<2$, we prove that, almost surely,
all unordered eigenfunction-product measures in every finite simple spectral
block are linearly independent. This strengthens the corresponding diagonal
independence statement and implies surjectivity of the joint eigenvalue and
eigenframe response map. We use this transversality to construct real-analytic
local deformations that preserve a prescribed finite block of eigenvalues
while varying its eigenframe. We then identify the least Cameron--Martin energy
required to realize a given infinitesimal angular response, expressing it
through a positive Schur complement of the response operator. Finally, on
finite-dimensional Gaussian slices, we derive the conditional law associated
with the local isospectral leaf and compare the induced slice metric with the
full Cameron--Martin minimum. The results provide a local geometric description
of spectral flexibility for Liouville Brownian motion.
\end{abstract}


\begingroup
\small
\tableofcontents
\endgroup

\section{Introduction}

The response of a spectral block to a change of medium includes both
eigenvalue motion and the motion of its eigenfunctions.  Prescribing this
response requires understanding which variations can be produced
independently.  A related variational question asks for the least energy
needed to turn the eigenframe while keeping the selected eigenvalues
fixed.  For a singular random speed measure, these questions connect
spectral perturbation theory with the geometry of the underlying Gaussian
field.

Liouville Brownian motion (LBM), constructed independently by Berestycki and
by Garban--Rhodes--Vargas, is the diffusion associated with Liouville quantum
gravity area \cite{BerestyckiDiffusion,GarbanRhodesVargas}. Its analytic
foundations were developed through the associated Dirichlet form, resolvent,
and heat kernel \cite{GarbanRhodesVargasHeatKernel,AndresKajino}, leading in
particular to spectral dimension two
\cite{RhodesVargasSpectral,AndresKajino}. Recent work also studies processes
time-changed by Gaussian multiplicative chaos, including a scaling-limit result
for LBM in the $L^2$ regime \cite{OoiGMC}. Berestycki and Wong established the
Weyl law for the LBM spectrum \cite{BerestyckiWong}; see
\cite{BerestyckiSpectral} for the broader emerging spectral picture. A
response calculus under Cameron--Martin perturbations then gave spectral
simplicity and absolute continuity of finite eigenvalue vectors
\cite{CaiResponse}.  These results describe global spectral growth and the
first-order response of eigenvalues.  We study instead the joint response
geometry of a finite spectral block, including the motion of its eigenframe.
Product response measures determine both the eigenvalue motion and the
rotation of the eigenframe.  Fixing the eigenvalues gives a minimum-energy
metric on the angular response.

Independence of eigenfunction products and local spectral flexibility
have substantial antecedents.  Privat and Sigalotti proved generic linear
independence of squared Dirichlet eigenfunctions under domain perturbations
\cite{PrivatSigalotti,PrivatSigalottiErratum}.  Pol\'a\v cik and Valdebenito
proved generic algebraic independence of radial bound-state eigenfunctions
for radial Schr\"odinger potentials \cite{PolacikValdebenito}.  For a
finite basis with independent products, realizing arbitrary matrix
elements by a local potential is a classical construction
\cite{Harriman}.  Product relations also govern symmetric-matrix
submersions in Strong Arnold transversality; see
\cite{GreilhuberKepplinger} and the earlier works discussed there.

Our main transversality assertion concerns all unordered product measures
in an arbitrary finite simple spectral block of LBM, almost surely for
each fixed $0<\gamma<2$.  The perturbation changes a singular speed measure,
and the proof separates the singular and absolutely continuous terms in a
distributional product equation.  This probability-one statement for the GMC speed measure extends the
diagonal square response used in \cite{CaiResponse}.  Perturbation theory,
constrained Hilbert-space minimization, and Gaussian disintegration then
yield local deformations, a response metric, and its conditional Gaussian
realization.

\subsection{The model and product response}

Let $D\subset\R^2$ be a bounded simply connected domain, let $h$ be a
zero-boundary Gaussian free field (GFF) on $D$ with law $\Prob$, and fix
$0<\gamma<2$.  Write $M_h$ for its subcritical Gaussian multiplicative chaos
measure.  The Cameron--Martin space is
\[
 \Hc=H_0^1(D),
 \qquad
 \inner{u}{v}_{\Hc}
 =\E(u,v)
 =\frac1{2\pi}\int_D\nabla u\cdot\nabla v\,dx.
\]
Let $A_h$ be the positive self-adjoint Dirichlet operator of Liouville
Brownian motion in $L^2(M_h)$, so that the Markov generator is $-A_h$.
The spectral response calculus of \cite{CaiResponse} gives, almost surely,
a simple spectrum and real $L^2(M_h)$-orthonormal eigenfunctions:
\[
 0<\lambda_1<\lambda_2<\cdots\uparrow\infty,
 \qquad
 A_h\phi_n=\lambda_n\phi_n.
\]
We use their bounded continuous representatives.  The vectors
$e_n=\phi_n/\sqrt{\lambda_n}$ form an orthonormal basis of $\Hc$.

For a functional $F$ of the field and a smooth Cameron--Martin direction
$f\in C_c^\infty(D)$, write
\[
 D_fF(h)=\left.\frac{d}{dt}F(h+tf)\right|_{t=0}
\]
whenever this derivative exists.  In particular, \cite{CaiResponse} gives
\begin{equation}
 \label{eq:intro-eigenvalue-response}
 D_f\log\lambda_n=-\gamma\int_Df\phi_n^2\,dM_h.
\end{equation}
Fix a nonempty finite block $I=\{n_1<\cdots<n_r\}$ of spectral labels.
Let $\Sym(r)$ denote the real symmetric $r\times r$ matrices, and define
\[
 \mathcal R_If
 =\left(\int_Df\phi_{n_a}\phi_{n_b}\,dM_h\right)_{a,b=1}^r
 \in\Sym(r).
\]
Spectral approximation of eigenfunction products is studied in
\cite{LuSoggeSteinerberger}.  Here the signed measures $\phi_i\phi_jM_h$ are the
coordinates of the field response.  In the next theorem, $a_{ij}$ are
arbitrary real coefficients and the first equality is an equality of signed
measures.

\begin{theorem}
\label{thm:intro-transversality}
Almost surely, for every finite block $I$,
\begin{equation}
 \label{eq:intro-product-independence}
 \sum_{\substack{i,j\in I\\i\le j}}a_{ij}\phi_i\phi_jM_h=0
 \quad\Longrightarrow\quad
 a_{ij}=0\quad(i\le j),
\end{equation}
and
\begin{equation}
 \label{eq:intro-response-surjective}
 \mathcal R_I(C_c^\infty(D))=\Sym(r).
\end{equation}
\end{theorem}

The proof separates the parts of the distributional product equation that
are singular and absolutely continuous with respect to Lebesgue measure.
Joint eigenvalue densities
from \cite{CaiResponse} exclude coincidences among distinct unordered pair
sums $\lambda_i+\lambda_j$.  Weak unique continuation rules out identically
vanishing products; Appendix~\ref{app:ucp} proves the required measure-valued
version using Sawyer's estimate \cite{SawyerUCP}.

\subsection{The isospectral response metric}

Write $\mathrm{Id}_{\Hc}$ for the identity operator on $\Hc$,
$\mathrm{Id}_r$ for the identity matrix, and
$\varepsilon_1,\ldots,\varepsilon_r$ for the standard basis of $\R^r$.
The rotation group and its tangent space at the identity are
\[
 \SO(r)=\{Q\in\R^{r\times r}:Q^TQ=\mathrm{Id}_r,\ \det Q=1\},
 \qquad
 \so(r)=\{A\in\R^{r\times r}:A^T=-A\}.
\]
We equip $\so(r)$ with
$\inner{A}{B}_{\so(r)}=\tfrac12\Tr(A^TB)$.
Thus $Q\in\SO(r)$ is a finite rotation, whereas $A\in\so(r)$ is an
infinitesimal rotation.

Along a smooth field perturbation, choose the signs of the energy-normalized
eigenfunctions continuously.  For $m\ne n$ their angular response is
\begin{equation}
 \label{eq:intro-connection}
 \omega_{mn}(f)
 :=\inner{e_m}{D_fe_n}_{\Hc}
 =\gamma\frac{\sqrt{\lambda_m\lambda_n}}{\lambda_m-\lambda_n}
   \int_Df\phi_m\phi_n\,dM_h,
\end{equation}
and $\omega_{nn}(f)=0$.  Set
\[
 \mathcal D_I(f)=(D_f\log\lambda_{n_a})_{a=1}^r,
 \qquad
 \Omega_I(f)=(\omega_{n_an_b}(f))_{a,b=1}^r,
 \qquad
 \mathcal C_I=(\mathcal D_I,\Omega_I).
\]
Finite Green energy of the product measures extends these maps
continuously to $\Hc$; see Section~\ref{sec:prelim}.  Theorem~\ref{thm:intro-transversality}
and simplicity of the spectrum then give
\begin{equation}
 \label{eq:intro-C-surjective}
 \mathcal C_I(C_c^\infty(D))
 =\mathcal C_I(\Hc)=\R^r\oplus\so(r).
\end{equation}
Let $P_I^0$ be the orthogonal projection of $\Hc$ onto
$\ker\mathcal D_I$.  The angular Gram operator is
\begin{equation}
 \label{eq:intro-angular-gram}
 K_I=\Omega_IP_I^0\Omega_I^*.
\end{equation}
Here adjoints use the specified Hilbert inner products.
Surjectivity in \eqref{eq:intro-C-surjective} makes
$\Omega_I|_{\ker\mathcal D_I}$ onto, hence $K_I$ is positive definite.
For $A\in\so(r)$ define
\begin{equation}
 \label{eq:intro-minimal-lift}
 H_I(A)=P_I^0\Omega_I^*K_I^{-1}A.
\end{equation}

The full response Gram operator has the block form
\[
 G_I=\mathcal C_I\mathcal C_I^*
 =\begin{pmatrix}A_I&B_I\\B_I^*&C_I\end{pmatrix},
\]
and its Schur complement is
\begin{equation}
 \label{eq:intro-schur-complement}
 K_I=C_I-B_I^*A_I^{-1}B_I.
\end{equation}

To describe finite perturbations, put $N=r(r+1)/2$ and choose
$\Hc$-orthonormal directions $f_1,\ldots,f_N\in C_c^\infty(D)$.
For $t=(t_1,\ldots,t_N)\in\R^N$ near $0$, set
\[
 h_t=h+\sum_{\alpha=1}^Nt_\alpha f_\alpha,
 \qquad
 E_I(t)\varepsilon_a=e_{n_a}(t).
\]
The selected eigenvalues $\lambda_{n_a}(t)$ and energy frame $E_I(t)$ are
real analytic near $0$.  The overlap matrix, its orthogonal polar factor,
and the corresponding angular coordinate are
\[
 M(t)=E_I(0)^*E_I(t),
 \qquad
 Q(t)=M(t)\bigl(M(t)^TM(t)\bigr)^{-1/2}\in\SO(r),
 \qquad
 \xi(t)=\log Q(t)\in\so(r),
\]
where $\log$ is the matrix logarithm near $\mathrm{Id}_r$.
Define the local response map
\begin{equation}
 \label{eq:intro-response-chart}
 \Phi_I(t)
 =\left((\log\lambda_{n_a}(t)-\log\lambda_{n_a}(0))_{a=1}^r,\xi(t)\right).
\end{equation}

\begin{theorem}
\label{thm:intro-connection-chart}
Almost surely, for every finite block $I$ and $A\in\so(r)$,
\begin{equation}
 \label{eq:intro-lift-identities}
 \mathcal D_I(H_I(A))=0,
 \qquad
 \Omega_I(H_I(A))=A,
\end{equation}
and $H_I(A)$ uniquely attains the minimum
\begin{equation}
 \label{eq:intro-lift-cost}
 \min_{\substack{f\in\Hc\\\mathcal C_I(f)=(0,A)}}\|f\|_{\Hc}^2
 =\|H_I(A)\|_{\Hc}^2
 =\inner{K_I^{-1}A}{A}_{\so(r)}.
\end{equation}
The directions $f_1,\ldots,f_N$ can be chosen so that $\Phi_I$ is a
real-analytic local diffeomorphism at $0$.  For every sufficiently small
$A\in\so(r)$,
\begin{equation}
 \label{eq:intro-isospectral-chart}
 \lambda_{n_a}\bigl(\Phi_I^{-1}(0,A)\bigr)=\lambda_{n_a}(0),
 \qquad 1\le a\le r.
\end{equation}
\end{theorem}

The rest of the paper is organized as follows.
Section~\ref{sec:prelim} records the spectral, Green-potential, and
Gaussian-slice inputs.  Section~\ref{sec:products} proves
Theorem~\ref{thm:intro-transversality}.
Section~\ref{sec:connection} constructs the response metric and local
charts of Theorem~\ref{thm:intro-connection-chart}.
Section~\ref{sec:applications} applies the metric to one eigenline and
to the conditional Gaussian law on a finite-slice isospectral leaf.
Appendix~\ref{app:ucp} supplies the measure-valued weak
unique-continuation argument.

\section{Preliminaries}
\label{sec:prelim}

We retain the notation of Section~1 and record the inputs from
\cite{CaiResponse} used below.

\subsection{Spectral and measure inputs}

Write $\widetilde u$ for the quasi-continuous representative of $u\in\Hc$.
The trace form defining $A_h$ has domain
\cite[Proposition~6.7]{CaiResponse}
\begin{equation}
 \label{eq:trace-form}
 \mathcal F_h
 =\{u\in\Hc:\widetilde u\in L^2(M_h)\}.
\end{equation}
The measure $M_h$ charges no set of zero Sobolev capacity
\cite[Lemma~6.6]{CaiResponse}; all integrals of Sobolev functions against
$M_h$ use these representatives.  The background operator is
$\Lzero=-(2\pi)^{-1}\Delta_D$, with Dirichlet Green kernel $\Green$.

We work on a common probability-one event for the following pathwise
properties.  Write $\delta_{nm}$ for the Kronecker delta.
Theorem~1.8 and Lemma~6.27 of \cite{CaiResponse} give
\begin{equation}
 \label{eq:response-input-spectrum}
 0<\lambda_1<\lambda_2<\cdots\uparrow\infty,
 \qquad
 A_h\phi_n=\lambda_n\phi_n,
 \qquad
 \int_D\phi_n\phi_m\,dM_h=\delta_{nm},
\end{equation}
with real bounded continuous representatives $\phi_n$.
Every finite ordered eigenvalue vector has a density with
respect to Lebesgue measure \cite[Theorem~1.8]{CaiResponse}.

The measure properties \cite[Lemmas~6.6, 6.16, and~6.28]{CaiResponse} are
\begin{equation}
 \label{eq:measure-structure-prelim}
 0<M_h(D)<\infty,
 \qquad
 M_h\perp dx,
 \qquad
 \supp M_h=D,
\end{equation}
where support is taken relative to $D$.
The Frostman estimate \cite[Lemma~6.10]{CaiResponse} provides a deterministic
$\alpha=\alpha(\gamma)>0$ and an almost surely finite random constant
$C_{\rm Fr}$ such that
\begin{equation}
 \label{eq:frostman-prelim}
 M_h(D\cap B(x,r))\le C_{\rm Fr}r^\alpha,
 \qquad x\in D,\quad 0<r\le1.
\end{equation}
Write $\log^+r=\max\{\log r,0\}$ for $r>0$.
By \cite[proof of Lemma~6.11]{CaiResponse}, the Frostman bound gives
the logarithmic Kato estimate: for $0<r\le1$,
\begin{equation}
 \label{eq:kato-gmc}
 \begin{split}
 &\sup_{x\in D}\int_{D\cap B(x,r)}
 \left(1+\log^+\frac1{|x-y|}\right)M_h(dy)\\
 &\qquad\le C_{\rm Fr}r^\alpha
 \left(1+\log\frac1r+\frac1\alpha\right)
 \xrightarrow[r\downarrow0]{}0.
 \end{split}
\end{equation}

\subsection{Green potentials}

For a finite signed measure $\sigma$ of finite Green energy, meaning
$I_D(|\sigma|)<\infty$, write
\[
 I_D(\sigma)=\iint_{D\times D}\Green(x,y)\,\sigma(dx)\sigma(dy).
\]
By \cite[Lemma~6.14]{CaiResponse}, the Green potential $U_D\sigma\in\Hc$
is uniquely characterized by
\begin{equation}
 \label{eq:green-potential-duality}
 \inner{U_D\sigma}{f}_{\Hc}
 =\int_D\widetilde f\,d\sigma,
 \qquad f\in\Hc,
\end{equation}
and
\begin{equation}
 \label{eq:green-energy-isometry}
 \|U_D\sigma\|_{\Hc}^2=I_D(\sigma).
\end{equation}
For all spectral labels $i,j$, \cite[Proposition~A.1]{CaiResponse} gives
\begin{equation}
 \label{eq:pair-finite-green-energy}
 I_D(|\phi_i\phi_j|M_h)<\infty.
\end{equation}
We therefore write
\begin{equation}
 \label{eq:Uij-def-prelim}
 U_{ij}=U_D(\phi_i\phi_jM_h),
 \qquad
 \inner{U_{ij}}{f}_{\Hc}
 =\int_D\widetilde f\,\phi_i\phi_j\,dM_h.
\end{equation}

\subsection{Finite Gaussian slices}
\label{subsec:prelim-slices}

Fix deterministic $\Hc$-orthonormal functions
$f_1,\ldots,f_N\in C_c^\infty(D)$ and set
\[
 V=\spanop\{f_1,\ldots,f_N\},
 \qquad
 v_t=\sum_{a=1}^Nt_af_a,\quad t\in\R^N.
\]
Let $W$ be the isonormal Gaussian process over $\Hc$ representing $h$,
and define
\begin{equation}
 \label{eq:prelim-slice-coordinates}
 t_a=W(f_a),
 \qquad
 h^\perp=h-v_t.
\end{equation}
Write $\Prob_V^\perp$ for the law of $h^\perp$ and
$\gamma_N(dt)=(2\pi)^{-N/2}e^{-|t|^2/2}dt$ for standard Gaussian measure.
By \cite[Lemma~2.3]{CaiResponse}, $t$ has law $\gamma_N$ and is independent
of $h^\perp$.  For every bounded measurable $F$,
\begin{equation}
 \label{eq:prelim-gaussian-disintegration}
 \Exp F(h)
 =\int\!\int_{\R^N}
 F(h^\perp+v_t)\,\gamma_N(dt)\,\Prob_V^\perp(dh^\perp).
\end{equation}

For $\Prob_V^\perp$-almost every $h^\perp$,
\cite[Lemma~2.4]{CaiResponse} supplies a family $(M_t)_{t\in\R^N}$ with
\begin{equation}
 \label{eq:prelim-coherent-slice}
 dM_t=e^{\gamma(v_t-v_s)}\,dM_s,
 \qquad s,t\in\R^N,
\end{equation}
agreeing with the canonical GMC of $h^\perp+v_t$ for
$\gamma_N$-almost every $t$.  On compact parameter sets, the relative
densities and their inverses are uniformly bounded; their square roots
are smooth multipliers of $\Hc$.

Fix such a family, a parameter $t_0$, and a finite block $I$ whose
eigenvalues at $t_0$ are simple.  Write $\lambda_i(t)$ and $\phi_i(t)$
for the eigenpairs of the trace form associated with $M_t$, and use
the normalization $e_i(t)=\phi_i(t)/\sqrt{\lambda_i(t)}$ from Section~1.
The simple-eigenvalue case of \cite[Proposition~2.5]{CaiResponse},
after undoing its unitary transport, gives the following statement.

\begin{proposition}
\label{prop:prelim-analytic-slice}
On a neighborhood $U$ of $t_0$, the eigenbranches can be chosen so that
\begin{equation}
 \label{eq:prelim-analytic-energy-frame}
 t\longmapsto\lambda_i(t)\in(0,\infty),
 \qquad
 t\longmapsto e_i(t)\in\Hc,
 \qquad i\in I,
\end{equation}
are real analytic, and
\begin{equation}
 \label{eq:prelim-slice-eigenrelations}
 \inner{e_i(t)}{e_j(t)}_{\Hc}
 =\int_D\phi_i(t)\phi_j(t)\,dM_t
 =\delta_{ij},
 \qquad i,j\in I.
\end{equation}
\end{proposition}

\Needspace{12\baselineskip}
\section{Product response and finite-block transversality}
\label{sec:products}

\subsection{Product identity and nonvanishing}

For $u,v\in\Hc$, write
\[
 \GammaE(u,v)=\frac1{2\pi}\nabla u\cdot\nabla v\in L^1(D,dx).
\]

\begin{lemma}
\label{lem:product-identity}
Almost surely, for every $i,j\ge1$,
\begin{equation}
 \label{eq:product-identity}
 \Lzero(\phi_i\phi_j)
 =(\lambda_i+\lambda_j)\phi_i\phi_jM_h
  -2\GammaE(\phi_i,\phi_j)\,dx
 \qquad\text{in }\mathcal D'(D).
\end{equation}
\end{lemma}

\begin{proof}
Boundedness of the eigenfunctions and the Sobolev product rule give
\[
 \phi_i\phi_j\in\Hc,
 \qquad
 \nabla(\phi_i\phi_j)
 =\phi_i\nabla\phi_j+\phi_j\nabla\phi_i\in L^2(D;\R^2).
\]
For $\zeta\in C_c^\infty(D)$ and $\ell\in\{i,j\}$,
\[
 \zeta\phi_\ell\in\Hc,
 \qquad
 \|\zeta\phi_\ell\|_{L^2(M_h)}
 \le\|\zeta\|_\infty\|\phi_\ell\|_{L^2(M_h)}
 =\|\zeta\|_\infty.
\]
Thus $\zeta\phi_i,\zeta\phi_j\in\mathcal F_h$, and the weak eigenvalue
equations yield
\begin{align*}
 \langle\Lzero(\phi_i\phi_j),\zeta\rangle
 &=\E(\phi_i\phi_j,\zeta)\\
 &=\E(\phi_i,\zeta\phi_j)+\E(\phi_j,\zeta\phi_i)
   -2\int_D\zeta\GammaE(\phi_i,\phi_j)\,dx\\
 &=(\lambda_i+\lambda_j)\int_D\zeta\phi_i\phi_j\,dM_h
   -2\int_D\zeta\GammaE(\phi_i,\phi_j)\,dx.
\end{align*}
This proves \eqref{eq:product-identity}.
\end{proof}

\begin{lemma}
\label{lem:no-open-zero}
Almost surely, each $\phi_i$ has zero set with empty interior, and
\begin{equation}
 \label{eq:product-not-zero}
 \phi_i\phi_j\not\equiv0,
 \qquad i,j\ge1.
\end{equation}
\end{lemma}

\begin{proof}
For every $i\ge1$, Section~\ref{sec:prelim} gives $\lambda_i>0$ and
\[
 \phi_i\in\Hc\cap C(D),
 \qquad
 -\Delta\phi_i=2\pi\lambda_i\phi_iM_h,
 \qquad
 \int_D\phi_i^2\,dM_h=1,
\]
where the equation holds in distributions.  Together with
\eqref{eq:kato-gmc}, these verify the hypotheses of
Corollary~\ref{cor:lbm-eigenfunction-ucp}.  Each set $\{\phi_i\ne0\}$ is
therefore open and dense in $D$, so
\[
 \{\phi_i\ne0\}\cap\{\phi_j\ne0\}\ne\varnothing
 \quad\Longrightarrow\quad
 \phi_i\phi_j\not\equiv0.\qedhere
\]
\end{proof}

\subsection{Final proof of Theorem~\ref{thm:intro-transversality}}

Put
\[
 \mathcal P=\{(i,j)\in\mathbb N^2:i\le j\}.
\]

For distinct $(i,j),(k,\ell)\in\mathcal P$, put $J=\{i,j,k,\ell\}$
and $c_n=\delta_{ni}+\delta_{nj}-\delta_{nk}-\delta_{n\ell}$ for $n\in J$.
Then $c\ne0$. By \cite[Theorem~1.8]{CaiResponse}, the vector
$(\lambda_n)_{n\in J}$ has a density with respect to Lebesgue measure.
The hyperplane $\{x\in\R^J:\sum_{n\in J}c_nx_n=0\}$ has Lebesgue
measure zero, so
\[
 \Prob\{\lambda_i+\lambda_j=\lambda_k+\lambda_\ell\}=0.
\]
Taking the countable union over distinct pairs shows that, almost surely,
\begin{equation}
 \label{eq:pair-sum-nonresonance}
 \lambda_i+\lambda_j\ne\lambda_k+\lambda_\ell
 \qquad\bigl((i,j)\ne(k,\ell)\text{ in }\mathcal P\bigr).
\end{equation}

Fix a nonempty finite block $I=\{n_1<\cdots<n_r\}$ and put
\[
 \mathcal P_I=\mathcal P\cap I^2,
 \qquad
 N=|\mathcal P_I|=\frac{r(r+1)}2.
\]
We order $\mathcal P_I$ lexicographically.  For $\alpha=(i,j)\in\mathcal P_I$,
set
\begin{equation}
 \label{eq:pair-sum-and-product-indexed}
 \sigma_\alpha=\lambda_i+\lambda_j,
 \qquad
 p_\alpha=\phi_i\phi_j.
\end{equation}

Work on the common probability-one event of Section~\ref{sec:prelim} and
Lemma~\ref{lem:no-open-zero} and \eqref{eq:pair-sum-nonresonance}.  For real coefficients
$c=(c_\alpha)_{\alpha\in\mathcal P_I}$, suppose
\begin{equation}
 \label{eq:product-relation0}
 \sum_{\alpha\in\mathcal P_I}c_\alpha p_\alpha M_h=0.
\end{equation}
For $k\ge0$, define
\begin{equation}
 \label{eq:qk-product-def}
 z_k=\sum_{\alpha\in\mathcal P_I}c_\alpha\sigma_\alpha^kp_\alpha
 \in\Hc\cap C(D).
\end{equation}
Continuity and full support of $M_h$ give
\[
 z_0M_h=0
 \quad\Longrightarrow\quad
 z_0=0\quad\text{on }D.
\]
If $z_k=0$, Lemma~\ref{lem:product-identity} gives
\begin{equation}
 \label{eq:qk-product-split}
 0=\Lzero z_k
 =z_{k+1}M_h
  -2\sum_{\alpha=(i,j)\in\mathcal P_I}
     c_\alpha\sigma_\alpha^k\GammaE(\phi_i,\phi_j)\,dx.
\end{equation}
Both terms on the right are finite signed measures.  Since $M_h\perp dx$,
their singular and absolutely continuous parts vanish separately.  Thus
\[
 z_{k+1}M_h=0
 \quad\Longrightarrow\quad
 z_{k+1}=0\quad\text{on }D.
\]
Induction yields
\begin{equation}
 \label{eq:vandermonde-relations}
 z_k=0,
 \qquad 0\le k<N.
\end{equation}

For every $x\in D$, these identities form the Vandermonde system
\[
 \left(\sigma_\alpha^k\right)_{\substack{0\le k<N\\\alpha\in\mathcal P_I}}
 \left(c_\alpha p_\alpha(x)\right)_{\alpha\in\mathcal P_I}=0.
\]
By \eqref{eq:pair-sum-nonresonance},
\[
 \det\left(\sigma_\alpha^k\right)_{\substack{0\le k<N\\\alpha\in\mathcal P_I}}
 =\prod_{\substack{\alpha,\beta\in\mathcal P_I\\\alpha<\beta}}
   (\sigma_\beta-\sigma_\alpha)\ne0.
\]
Invertibility and Lemma~\ref{lem:no-open-zero} therefore give
\[
 c_\alpha p_\alpha\equiv0
 \quad\Longrightarrow\quad
 c_\alpha=0,
 \qquad \alpha\in\mathcal P_I.
\]
This proves \eqref{eq:intro-product-independence}.

For surjectivity, let $B=(b_{ab})_{a,b=1}^r\in\Sym(r)$ satisfy
$\Tr(B^T\mathcal R_If)=0$ for every $f\in C_c^\infty(D)$.  Then
\[
 0=\Tr(B^T\mathcal R_If)
 =\int_D f\left(
   \sum_{a=1}^r b_{aa}\phi_{n_a}^2
   +2\sum_{1\le a<b\le r}b_{ab}\phi_{n_a}\phi_{n_b}
  \right)dM_h.
\]
The signed measure in this integral is zero, so
\eqref{eq:intro-product-independence} implies
\[
 b_{aa}=0\quad(1\le a\le r),
 \qquad
 2b_{ab}=0\quad(1\le a<b\le r).
\]
Thus, with orthogonality taken in the trace inner product on $\Sym(r)$,
the preceding calculation gives, on the same probability-one event and
for every finite block $I$,
\begin{align*}
 \bigl(\mathcal R_I(C_c^\infty(D))\bigr)^\perp&=\{0\},\\
 \dim\mathcal R_I(C_c^\infty(D))
 &=\dim\Sym(r)=\frac{r(r+1)}2,\\
 \mathcal R_I(C_c^\infty(D))&=\Sym(r).
\end{align*}

\section{The spectral response connection}
\label{sec:connection}

Fix $I=\{n_1<\cdots<n_r\}$ on the probability-one event used in
Section~3.  All adjoints use the Hilbert inner products from Section~1.2.

\subsection{Response in the fixed energy space}

The weak eigenvalue equation gives
\begin{equation}
 \label{eq:energy-orthonormal}
 \E(e_n,e_m)=\delta_{nm}.
\end{equation}
With $\|u\|_{\mathcal F_h}^2=\E(u,u)+\|u\|_{L^2(M_h)}^2$, the spectral
theorem and $C_c^\infty(D)\subset\mathcal F_h$ give
\[
 \overline{\spanop\{\phi_n:n\ge1\}}^{\,\|\cdot\|_{\mathcal F_h}}
 =\mathcal F_h,
 \qquad
 \overline{\mathcal F_h}^{\,\|\cdot\|_{\Hc}}=\Hc.
\]
Hence $(e_n)_{n\ge1}$ is an orthonormal basis of $\Hc$.

For $u\in\mathcal F_h$, the variational characterization of $\lambda_1$
gives
\begin{equation}
 \label{eq:trace-ineq}
 \int_D\widetilde u^2\,dM_h
 \le\lambda_1^{-1}\E(u,u).
\end{equation}
By density, the map $u\mapsto\widetilde u$ on $\mathcal F_h$ extends to
a bounded operator $J_h:\Hc\to L^2(M_h)$.
For $u_k\in C_c^\infty(D)$ converging to $u$ in
$\Hc$, a subsequence satisfies
\[
 u_k\longrightarrow\widetilde u
 \quad\text{quasi-everywhere},
 \qquad
 J_hu_k\longrightarrow J_hu
 \quad M_h\text{-almost everywhere}.
\]
Since $M_h$ charges no set of zero capacity, $J_hu=\widetilde u$
$M_h$-almost everywhere.  Thus $\mathcal F_h=\Hc$ and
\eqref{eq:trace-ineq} holds for every $u\in\Hc$.

Define
\begin{equation}
 \label{eq:T-def}
 T_h=J_h^*J_h.
\end{equation}
For $v\in\Hc$,
\[
 \inner{T_he_n}{v}_{\Hc}
 =\int_D e_n\widetilde v\,dM_h
 =\lambda_n^{-1}\E(e_n,v),
\]
so
\begin{equation}
 \label{eq:T-spectrum}
 T_he_n=\lambda_n^{-1}e_n,
 \qquad
 \inner{T_hu}{v}_{\Hc}
 =\int_D\widetilde u\,\widetilde v\,dM_h.
\end{equation}
Because $\lambda_n^{-1}\to0$ and $(e_n)$ is complete, $T_h$ is positive,
self-adjoint, and compact.  The identity $J_h^*J_h=T_h$ also makes $J_h$
compact.

For $g\in L^\infty(M_h;\R)$, define the bounded operator $T_h[g]$ on
$\Hc$ by
\begin{equation}
 \label{eq:Tg-def}
 \inner{T_h[g]u}{v}_{\Hc}
 =\int_Dg\,\widetilde u\,\widetilde v\,dM_h.
\end{equation}
Cauchy--Schwarz and \eqref{eq:trace-ineq} give
\[
 |\inner{T_h[g]u}{v}_{\Hc}|
 \le\|g\|_\infty\|J_hu\|_{L^2(M_h)}\|J_hv\|_{L^2(M_h)}
 \le\lambda_1^{-1}\|g\|_\infty\|u\|_{\Hc}\|v\|_{\Hc}.
\]
Writing $\|\cdot\|_{\op}$ for the operator norm, we obtain
\begin{equation}
 \label{eq:Tg-bound}
 \|T_h[g]\|_{\op}\le\lambda_1^{-1}\|g\|_\infty.
\end{equation}

For $f\in C_c^\infty(D)$, consider the coherent perturbation
$dM_\tau=e^{\gamma\tau f}dM_h$ and its energy operator
\begin{equation}
 \label{eq:Tt-def}
 \inner{T_\tau u}{v}_{\Hc}
 =\int_D\widetilde u\,\widetilde v\,e^{\gamma\tau f}\,dM_h.
\end{equation}
The estimate
\[
 \sum_{k=0}^\infty\frac{|\gamma\tau|^k}{k!}
 \|T_h[f^k]\|_{\op}
 \le\lambda_1^{-1}e^{|\gamma\tau|\|f\|_\infty}
\]
gives, after complexification of $\Hc$, the locally operator-norm
convergent expansion
\begin{equation}
 \label{eq:T-series-prelim}
 T_\tau=\sum_{k=0}^\infty\frac{(\gamma\tau)^k}{k!}T_h[f^k].
\end{equation}
Thus $\tau\mapsto T_\tau$ is entire.  For real $\tau$,
\[
 T_\tau=J_h^*\bigl(e^{\gamma\tau f}\,\cdot\bigr)J_h
\]
is positive, self-adjoint, and compact, and differentiation at $\tau=0$
yields
\begin{equation}
 \label{eq:T-derivative}
 \inner{D_fT_hu}{v}_{\Hc}
 =\gamma\int_Df\,\widetilde u\,\widetilde v\,dM_h.
\end{equation}

Choose the normalized real analytic eigenbranches from
Proposition~\ref{prop:prelim-analytic-slice}.  Differentiating their
orthonormality gives
\begin{equation}
 \label{eq:prelim-skew-frame-derivative}
 \inner{e_m}{D_fe_n}_{\Hc}
 =-\inner{e_n}{D_fe_m}_{\Hc},
 \qquad
 \inner{e_n}{D_fe_n}_{\Hc}=0.
\end{equation}

\begin{lemma}
\label{lem:connection-formula}
For every $n$,
\begin{equation}
 \label{eq:log-response}
 D_f\log\lambda_n=-\gamma\int_Df\phi_n^2\,dM_h,
\end{equation}
and, for $m\ne n$,
\begin{equation}
 \label{eq:omega-formula}
 \omega_{mn}(f)
 =\gamma\frac{\sqrt{\lambda_m\lambda_n}}{\lambda_m-\lambda_n}
  \int_Df\phi_m\phi_n\,dM_h.
\end{equation}
\end{lemma}

\begin{proof}
Put $\mu_n=\lambda_n^{-1}$.  From \eqref{eq:T-derivative},
\begin{equation}
 \label{eq:connection-Tmatrix}
 \inner{e_m}{D_fT_he_n}_{\Hc}
 =\frac{\gamma}{\sqrt{\lambda_m\lambda_n}}
  \int_Df\phi_m\phi_n\,dM_h.
\end{equation}
Differentiating $T_he_n=\mu_ne_n$ gives
\[
 (D_fT_h)e_n+(T_h-\mu_n\mathrm{Id}_{\Hc})D_fe_n
 =(D_f\mu_n)e_n.
\]
Taking the $e_n$ component and using
$\inner{e_n}{D_fe_n}_{\Hc}=0$ gives
\[
        D_f\mu_n
        =\inner{e_n}{D_fT_he_n}_{\Hc}.
\]
Hence
\[
 D_f\log\lambda_n
 =-\frac{D_f\mu_n}{\mu_n}
 =-\gamma\int_Df\phi_n^2\,dM_h.
\]
For $m\ne n$,
\[
 (\mu_n-\mu_m)\inner{e_m}{D_fe_n}_{\Hc}
 =\inner{e_m}{D_fT_he_n}_{\Hc}.
\]
Since
\[
        \mu_n-\mu_m
        =\frac{\lambda_m-\lambda_n}{\lambda_m\lambda_n},
\]
\eqref{eq:omega-formula} follows from
\eqref{eq:connection-Tmatrix}.
\end{proof}

By the Green-potential duality \eqref{eq:Uij-def-prelim}, the covectors in
\eqref{eq:log-response}--\eqref{eq:omega-formula} extend continuously from
$C_c^\infty(D)$ to $\Hc$.  With the notation of Section~1.2, their
extensions satisfy, for $f\in\Hc$ and $1\le a,b\le r$,
\begin{align}
 \label{eq:D-map}
        \mathcal D_I(f)_a
        &=-\gamma\inner{U_{n_an_a}}{f}_{\Hc},\\
 \label{eq:Omega-map}
        \Omega_I(f)_{ab}
        &=\gamma\frac{\sqrt{\lambda_{n_a}\lambda_{n_b}}}
        {\lambda_{n_a}-\lambda_{n_b}}
        \inner{U_{n_an_b}}{f}_{\Hc},
        \qquad a\ne b,
\end{align}
with $\Omega_I(f)_{aa}=0$.  Thus
\begin{equation}
 \label{eq:C-map-body}
        \mathcal C_I=(\mathcal D_I,\Omega_I):
        \Hc\longrightarrow\R^r\oplus\so(r)
\end{equation}
is a bounded linear map.

For $f\in C_c^\infty(D)$, let $B=\mathcal R_I f\in\Sym(r)$.  Then
\begin{equation}
 \label{eq:R-to-C-coordinate-map}
 \mathcal D_I(f)_a=-\gamma B_{aa},
 \qquad
 \Omega_I(f)_{ab}
 =\gamma\frac{\sqrt{\lambda_{n_a}\lambda_{n_b}}}
 {\lambda_{n_a}-\lambda_{n_b}}B_{ab}
 \quad(a\ne b).
\end{equation}
The coordinate map in \eqref{eq:R-to-C-coordinate-map} is a linear
isomorphism
\[
        \Sym(r)\longrightarrow\R^r\oplus\so(r),
\]
because $\gamma>0$ and the block spectrum is simple.
Theorem~\ref{thm:intro-transversality} therefore gives
\begin{equation}
 \label{eq:C-surjective-body}
 \mathcal C_I(C_c^\infty(D))
 =\mathcal C_I(\Hc)
 =\R^r\oplus\so(r).
\end{equation}

\subsection{Minimum-energy lifts and the response metric}
\label{subsec:schur-response}

Recall $P_I^0$ and $K_I=\Omega_IP_I^0\Omega_I^*$ from
\eqref{eq:intro-angular-gram}.  Surjectivity in
\eqref{eq:C-surjective-body} implies
\begin{equation}
 \label{eq:restricted-angular-map}
 \Omega_I(\ker\mathcal D_I)=\so(r).
\end{equation}

\begin{lemma}
\label{lem:angular-gram-positive}
For $A\in\so(r)\setminus\{0\}$,
\[
 \inner{K_IA}{A}_{\so(r)}
 =\|P_I^0\Omega_I^*A\|_{\Hc}^2>0.
\]
\end{lemma}

\begin{proof}
The equality follows from $(P_I^0)^*= (P_I^0)^2=P_I^0$.
If $P_I^0\Omega_I^*A=0$, choose $f\in\ker\mathcal D_I$ with
$\Omega_I(f)=A$.  Then
\[
 0=\inner{P_I^0\Omega_I^*A}{f}_{\Hc}
 =\inner{A}{\Omega_I(f)}_{\so(r)}
 =\|A\|_{\so(r)}^2.
\]
\end{proof}

The restricted angular map
$T_I:=\Omega_I|_{\ker\mathcal D_I}$ is bounded and surjective by
\eqref{eq:restricted-angular-map}. For a bounded surjective operator
$T:\mathcal H_0\to\so(r)$ between Hilbert spaces, the unique
minimum-norm solution of $Tf=A$ is $T^*(TT^*)^{-1}A$. Here
$T_I^*=P_I^0\Omega_I^*$ and $T_IT_I^*=K_I$, so
\begin{equation}
 \label{eq:minimal-lift-cost-body}
 \begin{gathered}
 H_I(A)=P_I^0\Omega_I^*K_I^{-1}A,
 \\
 \min_{\substack{f\in\Hc\\\mathcal D_If=0,\ \Omega_If=A}}
 \|f\|_{\Hc}^2
 =\|H_I(A)\|_{\Hc}^2
 =\inner{K_I^{-1}A}{A}_{\so(r)}.
 \end{gathered}
\end{equation}
The minimizer is unique, and its responses satisfy
\[
 \mathcal D_I(H_I(A))=0,
 \qquad
 \Omega_I(H_I(A))=\Omega_IP_I^0\Omega_I^*K_I^{-1}A=A.
\]
Thus the lift $H_I$ in \eqref{eq:intro-minimal-lift} is well defined.
For $A,B\in\so(r)$, its response metric is
\begin{equation}
 \label{eq:response-metric}
        g_I^{\rm resp}(A,B)
        =\inner{K_I^{-1}A}{B}_{\so(r)}.
\end{equation}

In particular,
\begin{equation}
 \label{eq:minimal-lift-gram}
 \inner{H_I(A)}{H_I(B)}_{\Hc}
 =\inner{K_I^{-1}A}{B}_{\so(r)}
 =g_I^{\rm resp}(A,B).
\end{equation}

The full response Gram operator has blocks
\[
 A_I=\mathcal D_I\mathcal D_I^*,
 \qquad
 B_I=\mathcal D_I\Omega_I^*,
 \qquad
 C_I=\Omega_I\Omega_I^*.
\]
Thus, relative to $\R^r\oplus\so(r)$,
\begin{equation}
 \label{eq:full-response-gram-block}
 G_I=\mathcal C_I\mathcal C_I^*
 =\begin{pmatrix}A_I&B_I\\B_I^*&C_I\end{pmatrix}>0.
\end{equation}
Surjectivity in \eqref{eq:C-surjective-body} also gives $A_I>0$.
The orthogonal-projection formula associated with the same
Moore--Penrose inverse gives
\begin{equation}
 \label{eq:kernel-projection-formula}
 P_I^0=\mathrm{Id}_{\Hc}-\mathcal D_I^*A_I^{-1}\mathcal D_I.
\end{equation}
The Schur-complement identities \cite[Theorem~1.3.3 and Section~1.5]{Bhatia}
then read
\begin{align}
 \label{eq:schur-complement}
 K_I&=C_I-B_I^*A_I^{-1}B_I,\\
 \label{eq:schur-inverse-corner}
 (G_I^{-1})_{\rm ang,ang}&=K_I^{-1}.
\end{align}

For a constant signed permutation matrix $S$, changing the frame to
$E_I'=E_IS$ gives $P_I^{0\prime}=P_I^0$ and
\begin{align}
 \label{eq:gauge-connection}
 \Omega_I'(f)&=S^T\Omega_I(f)S,\\
 \label{eq:gauge-K}
 K_I'(S^TAS)&=S^T(K_IA)S,
 \qquad H_I'(S^TAS)=H_I(A),\\
 \label{eq:gauge-metric}
 g_I^{\rm resp\,'}(S^TAS,S^TBS)&=g_I^{\rm resp}(A,B).
\end{align}
Thus the cost is independent of the signs and ordering of the block
eigenfunctions.

\subsection{Local isospectral response charts}
\label{subsec:connection-chart}

Set $N=r+\dim\so(r)=r(r+1)/2$ and choose a deterministic sequence
$(\zeta_k)_{k\ge1}\subset C_c^\infty(D)$ dense in $\Hc$.
Applying Gram--Schmidt to its linearly independent $N$-tuples gives a
countable deterministic family of $N$-dimensional subspaces, each with a
smooth $\Hc$-orthonormal basis.  The
coherent-slice conclusions of Section~\ref{subsec:prelim-slices} hold
simultaneously for this family on a probability-one event.

By \eqref{eq:C-surjective-body}, there are $g_1,\ldots,g_N\in C_c^\infty(D)$
whose responses form a basis of $\R^r\oplus\so(r)$.  In orthonormal
coordinates on this space,
\[
 \det\bigl[\mathcal C_Ig_1\ \cdots\ \mathcal C_Ig_N\bigr]\ne0.
\]
By continuity, the determinant remains nonzero after replacing each
$g_\alpha$ by a sufficiently close $\zeta_{k_\alpha}$.  Gram--Schmidt preserves
their span, so the resulting directions $f_1,\ldots,f_N$ satisfy
\begin{equation}
 \label{eq:chart-orthonormal-directions}
        V=\spanop\{f_1,\ldots,f_N\},
        \qquad
        \mathcal C_I|_{V}:V\xrightarrow{\;\sim\;}
        \R^r\oplus\so(r),
\end{equation}
with $f_1,\ldots,f_N\in C_c^\infty(D)$ orthonormal in $\Hc$.
Fix this slice.

For $t\in\R^N$ near $0$, use the perturbation and energy frame
\[
 h_t=h+\sum_{\alpha=1}^Nt_\alpha f_\alpha,
 \qquad
 E_I(t)\varepsilon_a=e_{n_a}(t).
\]
Proposition~\ref{prop:prelim-analytic-slice} gives a neighborhood of $0$
where these eigenbranches are simple and real analytic, with
$E_I(t)^*E_I(t)=\mathrm{Id}_r$.
The overlap $M(t)=E_I(0)^*E_I(t)$ satisfies $M(0)=\mathrm{Id}_r$.
On a sufficiently small neighborhood, $\det M(t)>0$ and the maps
\begin{align}
 \label{eq:polar-frame}
 Q(t)&=M(t)\bigl(M(t)^TM(t)\bigr)^{-1/2}\in\SO(r),\\
 \label{eq:xi-frame}
 \xi(t)&=\log Q(t)\in\so(r)
\end{align}
are real analytic, where $\log$ is the matrix logarithm near
$\mathrm{Id}_r$.  Hence the response map
\begin{equation}
 \label{eq:Phi-chart}
 \Phi_I(t)
 =\left((\log\lambda_{n_a}(t)-\log\lambda_{n_a}(0))_{a=1}^r,\xi(t)\right)
\end{equation}
is real analytic and satisfies $\Phi_I(0)=0$.
For $s\in\R^N$, write $v_s=\sum_{\alpha=1}^Ns_\alpha f_\alpha\in V$.

\begin{lemma}
\label{lem:chart-derivative}
At the chart center,
\begin{equation}
 \label{eq:chart-derivative}
 D\Phi_I(0)s=\mathcal C_I(v_s).
\end{equation}
\end{lemma}

\begin{proof}
The eigenvalue coordinates give $\mathcal D_I(v_s)$ by
Lemma~\ref{lem:connection-formula}.  For the overlap matrix,
\begin{align*}
 DM(0)s&=E_I(0)^*D_{v_s}E_I(0)=\Omega_I(v_s),\\
 D(M^TM)(0)s&=\Omega_I(v_s)^T+\Omega_I(v_s)=0.
\end{align*}
Differentiating \eqref{eq:polar-frame} therefore yields
\[
 DQ(0)s=DM(0)s=\Omega_I(v_s).
\]
Since $D(\log)_{\mathrm{Id}_r}$ is the identity on $\so(r)$,
\[
 D\xi(0)s=DQ(0)s=\Omega_I(v_s).
\]
Combining the two coordinates proves \eqref{eq:chart-derivative}.
\end{proof}

The identification of $D\xi$ with $\Omega_I$ in this lemma is at the chart
center.  Away from that point, $\xi$ is defined relative to the fixed
reference frame $E_I(0)$, whereas $\Omega_I=E_I^*dE_I$ is the connection
of the moving spectral frame.

\subsection{Final proof of Theorem~\ref{thm:intro-connection-chart}}

The surjectivity \eqref{eq:C-surjective-body} and the minimum-norm
formula \eqref{eq:minimal-lift-cost-body} give
\eqref{eq:intro-C-surjective} and
\eqref{eq:intro-lift-identities}--\eqref{eq:intro-lift-cost}.

For the orthonormal slice in
\eqref{eq:chart-orthonormal-directions},
Lemma~\ref{lem:chart-derivative} identifies $D\Phi_I(0)$ with
$\mathcal C_I|_V$ through the isometry $s\mapsto v_s$.
It is therefore an isomorphism.  The analytic inverse function theorem
gives an open neighborhood $U\subset\R^N$ of $0$ such that
\[
 \Phi_I:U\xrightarrow{\;\sim\;}\Phi_I(U)
\]
is a real-analytic diffeomorphism onto an open neighborhood of $0$.
For every sufficiently small
$A\in\so(r)$,
\[
        t_A=\Phi_I^{-1}(0,A)
\]
satisfies
\[
        \Phi_I(t_A)=(0,A),
        \qquad
        \lambda_{n_a}(t_A)=\lambda_{n_a}(0)
        \quad(1\le a\le r).
\]
Taking the countable intersection over finite blocks gives the common
probability-one event in the theorem.

\section{Application of the response metric}
\label{sec:applications}

We apply Theorem~\ref{thm:intro-connection-chart} to the within-block
response of one eigenline and to Gaussian conditioning on a finite-slice
isospectral leaf.

\subsection{Minimum response cost for one eigenline}
\label{subsec:star-cost}

Assume $r\ge2$ and fix $n=n_\ell\in I$.
Let $(\delta_j)_{j\in I\setminus\{n\}}$ be the standard basis of
$\R^{I\setminus\{n\}}$.  Define
\begin{equation}
 \label{eq:star-projection}
 \pi_{n,I}:\so(r)\longrightarrow\R^{I\setminus\{n\}},
 \qquad
 \pi_{n,I}(A)_{n_b}=(A\varepsilon_\ell)_b
 \quad(b\ne\ell).
\end{equation}
The inner product on $\so(r)$ from Section~1.2 gives, for $b,c\ne\ell$,
\[
 \pi_{n,I}^*\delta_{n_b}
 =\varepsilon_b\varepsilon_\ell^T
  -\varepsilon_\ell\varepsilon_b^T,
 \qquad
 \inner{\pi_{n,I}^*\delta_{n_b}}{\pi_{n,I}^*\delta_{n_c}}_{\so(r)}
 =\delta_{bc}.
\]
Consequently, for $f\in\Hc$ and $j\in I\setminus\{n\}$,
\begin{equation}
 \label{eq:star-projection-identities}
 \pi_{n,I}\pi_{n,I}^*=\mathrm{Id}_{\R^{I\setminus\{n\}}},
 \qquad
 \pi_{n,I}(\Omega_I(f))_j=\omega_{jn}(f).
\end{equation}
The compressed isospectral Gram operator is
\begin{equation}
 \label{eq:star-isospectral-gram}
 K_{n,I}^\star=\pi_{n,I}K_I\pi_{n,I}^*.
\end{equation}
For $0\ne x\in\R^{I\setminus\{n\}}$, the isometry of $\pi_{n,I}^*$
and $K_I>0$ give
\[
 x^TK_{n,I}^\star x
 =\inner{K_I\pi_{n,I}^*x}{\pi_{n,I}^*x}_{\so(r)}>0.
\]
Compressing the Schur identity \eqref{eq:schur-complement} gives
\begin{equation}
 \label{eq:star-schur-complement}
 K_{n,I}^\star
 =\pi_{n,I}C_I\pi_{n,I}^*
  -\pi_{n,I}B_I^*A_I^{-1}B_I\pi_{n,I}^*.
\end{equation}

\begin{proposition}
\label{prop:single-eigenline-cost}
For $x\in\R^{I\setminus\{n\}}$,
\begin{equation}
 \label{eq:star-response-metric-variational}
 \min_{\substack{f\in\Hc\\\mathcal D_If=0\\
                 \pi_{n,I}\Omega_If=x}}\|f\|_{\Hc}^2
 =\min_{\pi_{n,I}A=x}g_I^{\rm resp}(A,A)
 =x^T(K_{n,I}^\star)^{-1}x.
\end{equation}
Moreover,
\begin{align}
 (K_{n,I}^\star)^{-1}
 &\le\pi_{n,I}K_I^{-1}\pi_{n,I}^*,
 \label{eq:star-inverse-comparison}\\
 \Tr\bigl((K_{n,I}^\star)^{-1}\bigr)&\le\Tr(K_I^{-1}).
 \label{eq:star-cost-comparison}
\end{align}
\end{proposition}

\begin{proof}
Apply \eqref{eq:minimal-lift-cost-body} to the bounded surjection
$\pi_{n,I}\Omega_I|_{\ker\mathcal D_I}$, whose Gram operator is
$K_{n,I}^\star$.  The minimizing angular response is
\[
 A_x=K_I\pi_{n,I}^*(K_{n,I}^\star)^{-1}x,
 \qquad
 \pi_{n,I}A_x=x,
 \qquad
 g_I^{\rm resp}(A_x,A_x)=x^T(K_{n,I}^\star)^{-1}x.
\]
Together with \eqref{eq:minimal-lift-gram}, this gives
\eqref{eq:star-response-metric-variational}.
The admissible choice $A=\pi_{n,I}^*x$ gives
\[
 x^T(K_{n,I}^\star)^{-1}x
 \le\inner{K_I^{-1}\pi_{n,I}^*x}{\pi_{n,I}^*x}_{\so(r)},
\]
which is \eqref{eq:star-inverse-comparison}.  Since
$\pi_{n,I}^*\pi_{n,I}$ is an orthogonal projection,
\begin{align*}
 \Tr\bigl((K_{n,I}^\star)^{-1}\bigr)
 &\le\Tr(\pi_{n,I}K_I^{-1}\pi_{n,I}^*)\\
 &=\Tr(K_I^{-1}\pi_{n,I}^*\pi_{n,I})
 \le\Tr(K_I^{-1}).
\end{align*}
\end{proof}

For the unit response $x=\delta_j$, the unique minimizing lift is
\begin{equation}
 \label{eq:star-minimal-lifts}
 V_{nj}
 =H_I\bigl(K_I\pi_{n,I}^*(K_{n,I}^\star)^{-1}\delta_j\bigr)
 =P_I^0\Omega_I^*\pi_{n,I}^*(K_{n,I}^\star)^{-1}\delta_j.
\end{equation}
In particular,
\begin{equation}
 \label{eq:star-lift-identities}
 \mathcal D_I(V_{nj})=0,
 \qquad
 \omega_{kn}(V_{nj})=\delta_{kj}
 \quad(k,j\in I\setminus\{n\}),
\end{equation}
and
\begin{align}
 \|V_{nj}\|_{\Hc}^2
 &=\delta_j^T(K_{n,I}^\star)^{-1}\delta_j,
 \label{eq:star-individual-cost}\\
 \sum_{j\in I\setminus\{n\}}\|V_{nj}\|_{\Hc}^2
 &=\Tr\bigl((K_{n,I}^\star)^{-1}\bigr).
 \label{eq:star-total-cost}
\end{align}
These constraints fix the block eigenvalue responses and the angular
response of $e_n$ within the block.

For $x,y\in\R^{I\setminus\{n\}}$, the resulting response metric is
\begin{equation}
 \label{eq:star-response-metric}
 g_{n,I}^{\star}(x,y)=x^T(K_{n,I}^\star)^{-1}y.
\end{equation}

\subsection{Conditional Gaussian law on an isospectral slice}
\label{subsec:gaussian-leaves}

Fix a finite block $I=\{n_1<\cdots<n_r\}$, put $N=r(r+1)/2$, and choose
deterministic $\Hc$-orthonormal directions
$f_1,\ldots,f_N\in C_c^\infty(D)$ spanning $V$.
Condition on $h^\perp$ as in \eqref{eq:prelim-gaussian-disintegration} and
use the absolute Gaussian coordinates
\[
 h_t=h^\perp+v_t,
 \qquad v_t=\sum_{\alpha=1}^Nt_\alpha f_\alpha,
 \qquad t\sim\gamma_N.
\]
Write $\mathcal C_I(t)$, $K_I(t)$, and $H_{I,t}$ for the response objects
at $h_t$.  At a point $t_0$ where the block is simple and
\[
 \mathcal C_I(t_0)|_V:
 V\xrightarrow{\;\sim\;}\R^r\oplus\so(r),
\]
use the angular coordinate
\eqref{eq:polar-frame}--\eqref{eq:xi-frame} with overlap matrix
$E_I(t_0)^*E_I(t)$, and put
\[
 Y(t)=\bigl(\log\lambda_{n_a}(t)-\log\lambda_{n_a}(t_0)\bigr)_{a=1}^r.
\]
The chart construction in Section~\ref{sec:connection}, centered at $t_0$,
gives a real-analytic diffeomorphism
\begin{equation}
 \label{eq:leaf-product-chart}
 \Phi_I=(Y,\xi):U\xrightarrow{\;\sim\;}
 B_{\rm rad}\times B_{\rm ang},
 \qquad \Phi_I(t_0)=(0,0).
\end{equation}
Here $B_{\rm rad}\subset\R^r$ and $B_{\rm ang}\subset\so(r)$ are
origin-centered balls whose closed product lies inside the original
chart image.  The inverse
\begin{equation}
 \label{eq:inverse-response-chart-leaf}
 t=\Theta(y,\xi)=\Phi_I^{-1}(y,\xi)
\end{equation}
extends real analytically to a neighborhood of
$\overline{B_{\rm rad}}\times\overline{B_{\rm ang}}$.
Choose the balls small enough that $\mathcal C_I(t)|_V$ remains invertible
on $\overline U$.

For $y\in B_{\rm rad}$ put
\[
 \Sigma_y=Y^{-1}(y)\cap U,
 \qquad T_t\Sigma_y=\ker DY(t)\quad(t\in\Sigma_y).
\]
Since $D\Phi_I(t)$ is invertible, $DY(t)$ has rank $r$ on $U$.
Fix an orthonormal basis of $\so(r)$ and use the corresponding Lebesgue
measure $d\xi$.  The parametrization $\xi\mapsto\Theta(y,\xi)$ of
$\Sigma_y$ has induced metric
\begin{equation}
 \label{eq:leaf-metric}
 g_y(\xi)=D_\xi\Theta(y,\xi)^TD_\xi\Theta(y,\xi).
\end{equation}
Define the normal Jacobian and induced volume by
\begin{align}
 J_Y(t)&=\det\bigl(DY(t)DY(t)^T\bigr)^{1/2}>0,
 \label{eq:normal-jacobian}\\
 d\operatorname{vol}_{g_y}(\xi)&=\sqrt{\det g_y(\xi)}\,d\xi.
 \label{eq:leaf-volume}
\end{align}
Write $d\operatorname{vol}_{\Sigma_y}$ for the induced Riemannian volume
on $\Sigma_y$, and set
\begin{equation}
 \label{eq:leaf-normalization}
 Z_y=\int_{B_{\rm ang}}
 \frac{e^{-|\Theta(y,\eta)|^2/2}}{J_Y(\Theta(y,\eta))}
 \,d\operatorname{vol}_{g_y}(\eta).
\end{equation}

\begin{proposition}
\label{prop:gaussian-leaf-application}
For $\Prob_V^\perp$-almost every $h^\perp$, on each such chart $Z_y$ is
smooth and strictly positive, and the conditional law of $\xi$ given
$Y=y$ under $\gamma_N(\,\cdot\mid U)$ admits the version
\begin{equation}
 \label{eq:intro-conditional-leaf-law}
 d\nu_y(\xi)
 =Z_y^{-1}\frac{e^{-|\Theta(y,\xi)|^2/2}}{J_Y(\Theta(y,\xi))}
 \,d\operatorname{vol}_{g_y}(\xi),
 \qquad y\in B_{\rm rad}.
\end{equation}
At the chart center, for $A\in\so(r)$,
\begin{equation}
 \label{eq:intro-leaf-response-cost}
 g_0(0)(A,A)
 =\min_{\substack{f\in V\\\mathcal C_I(t_0)f=(0,A)}}\|f\|_{\Hc}^2
 \ge\inner{K_I(t_0)^{-1}A}{A}_{\so(r)}.
\end{equation}
Equality holds if and only if $H_{I,t_0}(A)\in V$.
\end{proposition}

\begin{proof}
For $\Prob_V^\perp$-almost every $h^\perp$, the coherent family
\eqref{eq:prelim-coherent-slice} contains a canonical measure
$M_{t_*}=M_{h_{t_*}}$, for some $t_*\in\R^N$, satisfying the pathwise
inputs of Section~\ref{sec:prelim}.
Work on such a fiber.  Since
\[
 dM_t=e^{\gamma(v_t-v_{t_*})}\,dM_{t_*},
 \qquad v_t-v_{t_*}\in C_c^\infty(D),
\]
the finite Green-energy estimate of \cite[Lemma~A.2]{CaiResponse} applies
at every $M_t$.  Where the block spectrum is simple, the response maps
\eqref{eq:D-map}--\eqref{eq:Omega-map} therefore extend continuously to
$\Hc$; write them as $\mathcal D_I(t)$, $\Omega_I(t)$, and
$\mathcal C_I(t)=(\mathcal D_I(t),\Omega_I(t))$.

When $r=1$, the angular space is $\so(1)=\{0\}$ and each leaf is a point.
With unit zero-dimensional volume and $\det g_y=1$, the ordinary
one-dimensional change-of-variables formula gives
\[
 \begin{gathered}
 \abs{\det D\Theta(y,0)}=J_Y(\Theta(y,0))^{-1},\\
 Z_y=\frac{e^{-|\Theta(y,0)|^2/2}}{J_Y(\Theta(y,0))},
 \qquad \nu_y=\delta_0.
 \end{gathered}
\]
Both variational minima in \eqref{eq:intro-leaf-response-cost} are attained
at $f=0$ and equal zero.  This proves the proposition when $r=1$.
For the rest of the proof assume $r\ge2$.

Since $DY$ has rank $r$ on a neighborhood of $\overline U$, the local
coarea formula \cite[Section~2.3, Proposition~2]{DiaconisHolmesShahshahani}
gives
\begin{equation}
 \label{eq:coarea-jacobian-identity}
 \abs{\det D\Theta(y,\xi)}
 =\frac{\sqrt{\det g_y(\xi)}}{J_Y(\Theta(y,\xi))}.
\end{equation}
In particular, the normalization in \eqref{eq:leaf-normalization} is
\[
 Z_y=\int_{B_{\rm ang}}e^{-|\Theta(y,\eta)|^2/2}
       \abs{\det D\Theta(y,\eta)}\,d\eta.
\]
The integrand is smooth and strictly positive on the closed product,
so $Z_y$ is smooth and $0<Z_y<\infty$ for every $y\in B_{\rm rad}$.

A subscript $\#$ denotes pushforward of measures.  Change of variables
gives
\begin{equation}
 \label{eq:response-coordinate-density}
 d((\Phi_I)_\#(\gamma_N|_U))(y,\xi)
 =(2\pi)^{-N/2}e^{-|\Theta(y,\xi)|^2/2}
   \abs{\det D\Theta(y,\xi)}\,dy\,d\xi.
\end{equation}
Its $y$-marginal and disintegration are
\begin{align*}
 Y_\#(\gamma_N|_U)(dy)&=(2\pi)^{-N/2}Z_y\,dy,\\
 (\Phi_I)_\#(\gamma_N|_U)(dy,d\xi)
 &=(2\pi)^{-N/2}Z_y\,dy\,\nu_y(d\xi).
\end{align*}
Thus the conditional laws admit the version
\begin{align}
 d\nu_y(\xi)
 &=Z_y^{-1}e^{-|\Theta(y,\xi)|^2/2}
   \abs{\det D\Theta(y,\xi)}\,d\xi
 \label{eq:conditional-leaf-density}\\
 &=Z_y^{-1}\frac{e^{-|\Theta(y,\xi)|^2/2}}
                     {J_Y(\Theta(y,\xi))}
   \,d\operatorname{vol}_{g_y}(\xi),
 \qquad y\in B_{\rm rad}.
 \label{eq:leaf-density-volume}
\end{align}
Equivalently, on $\Sigma_y$,
\begin{equation}
 \label{eq:coarea-leaf-law}
 d\bigl((\Theta(y,\cdot))_\#\nu_y\bigr)(t)
 =Z_y^{-1}\frac{e^{-|t|^2/2}}{J_Y(t)}
   \,d\operatorname{vol}_{\Sigma_y}(t).
\end{equation}
These formulas specify a smooth version for every $y\in B_{\rm rad}$.
Normalization by $\gamma_N(U)$ leaves the conditional laws unchanged.
This proves \eqref{eq:intro-conditional-leaf-law}.

For $t\in U$, the restriction $\mathcal C_I(t)|_V$ is invertible, so
$\mathcal C_I(t)$ is surjective.  The full projection $P_I^0(t)$, Gram
operator $K_I(t)$, and lift $H_{I,t}$ at $M_t$ are as in
\eqref{eq:intro-angular-gram}--\eqref{eq:intro-minimal-lift}.
For $s\in\R^N$, the eigenvalue response formula gives
\[
 DY(t)s=\mathcal D_I(t)v_s.
\]
Thus the isometry $s\mapsto v_s$ identifies the tangent space of the leaf
with
\[
 T_t^{\rm sl}
 =V\cap\ker\mathcal D_I(t)
 =\{v_s:s\in T_t\Sigma_{Y(t)}\}\subset\Hc.
\]
Let $P_{I,V}^0(t)$ be the orthogonal projection of $\Hc$ onto
$T_t^{\rm sl}$, and define the slice Gram operator
\begin{equation}
 \label{eq:slice-angular-gram}
        K_{I,V}(t)
        =\Omega_I(t)P_{I,V}^0(t)\Omega_I(t)^*.
\end{equation}
The restriction $\Omega_I(t)|_{T_t^{\rm sl}}$ is an isomorphism onto
$\so(r)$, so $K_{I,V}(t)>0$.  Since
$T_t^{\rm sl}\subset\ker\mathcal D_I(t)$, the orthogonal projections
satisfy
\[
        0\le P_{I,V}^0(t)\le P_I^0(t).
\]
Applying $\Omega_I(t)$ and its adjoint gives, in the sense of quadratic
forms,
\[
        0<K_{I,V}(t)\le K_I(t).
\]
For $B\in\so(r)$, put
\[
 v_B
 =P_{I,V}^0(t)\Omega_I(t)^*
   K_{I,V}(t)^{-1}B\in T_t^{\rm sl}.
\]
Then
\begin{align*}
 \Omega_I(t)v_B&=B,\\
 \|v_B\|_{\Hc}^2
 &=\inner{K_{I,V}(t)^{-1}B}
   {K_{I,V}(t)K_{I,V}(t)^{-1}B}_{\so(r)}
 =\inner{K_{I,V}(t)^{-1}B}{B}_{\so(r)}.
\end{align*}
Since $\Omega_I(t)|_{T_t^{\rm sl}}$ is injective, $v_B$ is the unique
admissible slice direction.  Hence
\begin{equation}
 \label{eq:slice-angular-cost}
 \min_{\substack{v\in T_t^{\rm sl}\\\Omega_I(t)v=B}}
   \|v\|_{\Hc}^2
 =\|v_B\|_{\Hc}^2
 =\inner{K_{I,V}(t)^{-1}B}{B}_{\so(r)}.
\end{equation}
The minimum-norm formula \eqref{eq:minimal-lift-cost-body}, applied to
the bounded surjective map $\mathcal C_I(t)$, gives
\begin{equation}
 \label{eq:full-angular-cost-pointwise}
 \min_{\substack{v\in\ker\mathcal D_I(t)\\\Omega_I(t)v=B}}
   \|v\|_{\Hc}^2
 =\inner{K_I(t)^{-1}B}{B}_{\so(r)}.
\end{equation}
The two lifts have the same response, so
\[
 \mathcal D_I(t)(v_B-H_{I,t}(B))=0,
 \qquad \Omega_I(t)(v_B-H_{I,t}(B))=0.
\]
Using $H_{I,t}(B)=P_I^0(t)\Omega_I(t)^*K_I(t)^{-1}B$, we obtain
\[
 \inner{H_{I,t}(B)}{v_B-H_{I,t}(B)}_{\Hc}
 =\inner{K_I(t)^{-1}B}{\Omega_I(t)(v_B-H_{I,t}(B))}_{\so(r)}=0.
\]
Consequently,
\[
 \|v_B\|_{\Hc}^2
 =\inner{K_I(t)^{-1}B}{B}_{\so(r)}
  +\|v_B-H_{I,t}(B)\|_{\Hc}^2.
\]
Thus,
\begin{equation}
 \label{eq:gram-loewner-comparison}
        0<K_{I,V}(t)\le K_I(t),
        \qquad
        K_I(t)^{-1}\le K_{I,V}(t)^{-1}.
\end{equation}
The two minimum costs agree exactly when $v_B=H_{I,t}(B)$, equivalently
$H_{I,t}(B)\in V$.

At $t_0=\Theta(0,0)$, fix $A\in\so(r)$ and write
$s_A=D_\xi\Theta(0,0)A\in\R^N$.
Differentiating $\Phi_I(\Theta(y,\xi))=(y,\xi)$ and using
Lemma~\ref{lem:chart-derivative} centered at $t_0$ gives
\[
 \mathcal D_I(t_0)v_{s_A}=DY(t_0)s_A=0,
 \qquad
 \Omega_I(t_0)v_{s_A}=D\xi(t_0)s_A=A.
\]
By uniqueness of the slice lift,
\[
 v_{s_A}
 =P_{I,V}^0(t_0)\Omega_I(t_0)^*K_{I,V}(t_0)^{-1}A.
\]
Hence, for $A,B\in\so(r)$,
\begin{align}
 g_0(0)(A,B)
 &=s_A\cdot s_B
 =\inner{v_{s_A}}{v_{s_B}}_{\Hc}\notag\\
 \label{eq:center-leaf-response-metric}
 &=\inner{K_{I,V}(t_0)^{-1}A}{B}_{\so(r)}.
\end{align}
Taking $B=A$ and using \eqref{eq:slice-angular-cost} and
\eqref{eq:gram-loewner-comparison} proves
\eqref{eq:intro-leaf-response-cost}, with equality exactly when
$H_{I,t_0}(A)\in V$.
\end{proof}

\appendix
\Needspace{12\baselineskip}
\section*{Appendix}

\section{Measure-valued weak unique continuation}
\label{app:ucp}

We prove the weak unique-continuation result used in
Lemma~\ref{lem:no-open-zero}, extending the absorption argument based on
Sawyer's logarithmic Taylor estimate \cite[(2.1)]{SawyerUCP} to measures.

Let $\Omega\subset\R^2$ be open, let $\mu$ be a locally finite positive
Radon measure on $\Omega$, and fix $a\in\R$.  Write $B_r=B(0,r)$.
For compact $K\Subset\Omega$ and sufficiently small $r>0$, put
\begin{equation}
 \label{eq:kato-modulus-app}
 \kappa_K(r)
 =\sup_{x\in K}\int_{B(x,r)}
   \left(1+\log^+\frac1{|x-y|}\right)\mu(dy).
\end{equation}
Throughout Appendix~\ref{app:ucp} we assume
\begin{equation}
 \label{eq:measure-kato-global}
 \kappa_K(r)\longrightarrow0
 \qquad(r\downarrow0),\qquad K\Subset\Omega.
\end{equation}
In particular, $\mu$ has no atoms: $\mu(\{x\})>0$ would make
$\kappa_{\{x\}}(r)=\infty$ for every sufficiently small $r>0$.

\subsection{Logarithmic Kato measures}
\label{subsec:kato-potential}

Fix a compact set $K_0\Subset U_1\Subset\Omega$, with $U_1$ open.
Choosing $\rho>0$ so that $\kappa_{K_0}(\rho)<\infty$ gives
\begin{equation}
 \label{eq:kato-uniform-log}
 \begin{split}
 &\sup_{x\in K_0}\int_{\overline U_1}
       \left(1+\log^+\frac1{|x-y|}\right)\mu(dy)\\
 &\qquad\le
 \kappa_{K_0}(\rho)
 +\left(1+\log^+\frac1\rho\right)\mu(\overline U_1)<\infty.
 \end{split}
\end{equation}
Let $b\in L^\infty(U_1,\mu)$ and $g\in L^2(U_1,dx)$, both real-valued
and extended by zero outside $U_1$.

\begin{lemma}
\label{lem:kato-stability-weights}
For all sufficiently small $0<r\le1$,
\begin{equation}
 \label{eq:kato-weighted-modulus}
 \sup_{x\in K_0}\int_{B(x,r)}
 \left(1+\log^+\frac1{|x-y|}\right)|b(y)|\,\mu(dy)
 \le \|b\|_\infty\kappa_{K_0}(r),
\end{equation}
and
\begin{equation}
 \label{eq:kato-L2-modulus}
 \sup_{x\in K_0}\int_{B(x,r)}
 \left(1+\log^+\frac1{|x-y|}\right)|g(y)|\,dy
 \le \sqrt{3\pi}\,\|g\|_2\,r(1+|\log r|).
\end{equation}
\end{lemma}

\begin{proof}
The first bound follows from $|b|\le\|b\|_\infty$ $\mu$-almost everywhere.
For the second, Cauchy--Schwarz gives
\begin{align*}
 &\int_{B(x,r)}
   \left(1+\log^+\frac1{|x-y|}\right)|g(y)|\,dy\\
 &\qquad\le
 \|g\|_2\left(2\pi\int_0^r t(1-\log t)^2\,dt\right)^{1/2},\\
 &2\pi\int_0^r t(1-\log t)^2\,dt
 =\pi r^2\left((1-\log r)^2+(1-\log r)+\tfrac12\right)
 \le3\pi r^2(1+|\log r|)^2.
\end{align*}
\end{proof}

Let $B\Subset\Omega$ be an open ball and let $u\in H^1(B)\cap C(B)$
be real-valued, with $-\Delta u=a u\mu$ in distributions on $B$.
For $\chi\in C_c^\infty(B)$, put
\[
 w=\chi u,\qquad
 g=-2\nabla\chi\cdot\nabla u-u\Delta\chi,
 \qquad \sigma=a\chi u\mu+g\,dx.
\]
Extend these functions and the measure by zero outside $B$.

\begin{lemma}
\label{lem:cutoff-laplacian-measure}
Then $\sigma$ is a finite signed Radon measure and
\begin{equation}
 \label{eq:cutoff-laplacian-general}
 w\in H^1(\R^2)\cap C_c(\R^2),
 \qquad -\Delta w=\sigma.
\end{equation}
Moreover,
\begin{equation}
 \label{eq:sigma-kato-newtonian}
 \sup_{x\in\R^2}\int_{B(x,\rho)}
 \left(1+\log^+\frac1{|x-y|}\right)|\sigma|(dy)
 \longrightarrow0
 \qquad(\rho\downarrow0).
\end{equation}
\end{lemma}

\begin{proof}
The Sobolev product rule gives \eqref{eq:cutoff-laplacian-general}, while
$u$ is bounded on $\supp\chi$ and
\[
 \|g\|_2
 \le2\|\nabla\chi\|_\infty\|\nabla u\|_{L^2(B)}
    +\|\Delta\chi\|_\infty\|u\|_{L^2(B)}<\infty.
\]
Both terms of $\sigma$ have compact support in $B$, so $\sigma$ is finite.
Choose a compact neighborhood $K_\chi\Subset B$ of $\supp\chi$.
For small $\rho$, only centers in $K_\chi$ contribute to the term
$|\chi u|\mu$, and Lemma~\ref{lem:kato-stability-weights} gives
\begin{align*}
 &\sup_{x\in\R^2}\int_{B(x,\rho)}
 \left(1+\log^+\frac1{|x-y|}\right)|\sigma|(dy)\\
 &\qquad\le
 |a|\,\|\chi u\|_\infty\kappa_{K_\chi}(\rho)
 +\sqrt{3\pi}\,\|g\|_2\,\rho(1+|\log\rho|)
 \longrightarrow0.
\end{align*}
\end{proof}

\subsection{The logarithmic Taylor kernel}

Let
\begin{equation}
 \label{eq:fundamental-solution}
 \Phi(x)=-\frac1{2\pi}\log|x|,
 \qquad -\Delta\Phi=\delta_0.
\end{equation}
For $y\ne0$ and an integer $N\ge0$, let $T_N(\cdot,y)$ be the Taylor
polynomial of degree $N$ at $x=0$ of $x\mapsto\Phi(x-y)$, and put
\begin{equation}
 \label{eq:taylor-remainder-kernel}
 R_N(x,y)=\Phi(x-y)-T_N(x,y).
\end{equation}
For $r>0$ define
\[
 L_r(x,y)=1+\log^+\frac{16r}{|x-y|},
 \qquad c_r=1+\log^+(16r).
\]
The dilation $(x,y)\mapsto(x/(16r),y/(16r))$ sends
$0<|x|,|y|<4r$ into Sawyer's range $0<|x|,|y|<1/4$.
The additive constant $(2\pi)^{-1}\log(16r)$ cancels from the Taylor
remainder.  Thus \cite[(2.1)]{SawyerUCP} gives an absolute constant $C>0$
such that, for $0<|x|,|y|<4r$ and $x\ne y$,
\begin{equation}
 \label{eq:sawyer-kernel-scaled}
 |R_N(x,y)|
 \le C\left(\frac{|x|}{|y|}\right)^{N+1}L_r(x,y),
 \qquad N\ge0.
\end{equation}
We keep this constant $C$ fixed below.

Let $w\in H^1(\R^2)\cap C_c(\R^2)$ be real-valued and suppose that
$\sigma=-\Delta w$ is a finite signed Radon measure satisfying
\eqref{eq:sigma-kato-newtonian}.  Assume $w=0$ on $B_\delta$ for some
$\delta>0$.

\begin{lemma}
\label{lem:newtonian-measure-representation}
For every $x\in\R^2$,
\begin{equation}
 \label{eq:newtonian-representation}
 w(x)=\int_{\R^2}\Phi(x-y)\,\sigma(dy),
\end{equation}
and
\begin{equation}
 \label{eq:taylor-subtracted-representation}
 w(x)=\int_{\R^2}R_N(x,y)\,\sigma(dy),
 \qquad N\ge0.
\end{equation}
\end{lemma}

\begin{proof}
Since $\supp\sigma\subset\supp w$ is compact, the splitting in
\eqref{eq:kato-uniform-log}, applied to $|\sigma|$, makes
$V=\Phi*\sigma$ finite everywhere.  For
$|z-x|<\rho<1/4$, splitting at $B(x,2\rho)$ also gives
\begin{align*}
 |V(z)-V(x)|
 &\le \frac1\pi\sup_{\xi\in\R^2}
   \int_{B(\xi,3\rho)}
   \left(1+\log^+\frac1{|\xi-y|}\right)|\sigma|(dy)\\
 &\quad+\frac{|z-x|}{2\pi\rho}|\sigma|(\R^2).
\end{align*}
First letting $z\to x$ and then $\rho\downarrow0$ proves continuity of $V$.
In distributions, $-\Delta V=\sigma=-\Delta w$.
Choose $\chi_0\in C_c^\infty(\R^2)$ equal to one near $\supp w$.
Then
\[
 \sigma(\R^2)
 =\langle-\Delta w,\chi_0\rangle
 =\int_{\R^2}\nabla w\cdot\nabla\chi_0\,dx=0.
\]
Consequently, for all sufficiently large $|x|$,
\[
 |V(x)|
 =\left|\int_{\R^2}[\Phi(x-y)-\Phi(x)]\,\sigma(dy)\right|
 \le\frac1{\pi|x|}\int_{\R^2}|y|\,|\sigma|(dy)
 \longrightarrow0.
\]
The continuous function $w-V$ is harmonic on $\R^2$ and tends to zero
at infinity; hence it vanishes, proving \eqref{eq:newtonian-representation}.

Since $w=0$ on $B_\delta$, we also have $|\sigma|(B_\delta)=0$.
The kernel is smooth in $x$ on $B_{\delta/2}$ for $y\in\supp\sigma$.
For every multi-index $\alpha$,
\[
 \int_{\R^2}(\partial^\alpha\Phi)(-y)\,\sigma(dy)
 =\partial^\alpha V(0)=\partial^\alpha w(0)=0.
\]
It follows that
\[
 \int_{\R^2}T_N(x,y)\,\sigma(dy)
 =\sum_{|\alpha|\le N}\frac{x^\alpha}{\alpha!}
   \int_{\R^2}(\partial^\alpha\Phi)(-y)\,\sigma(dy)=0.
\]
Subtracting this identity from \eqref{eq:newtonian-representation} proves
\eqref{eq:taylor-subtracted-representation}.
\end{proof}

Retain $w,\sigma,\delta$ as above and choose $r>\delta$ with
$\supp\sigma\subset B_{3r}$.

\begin{lemma}
\label{lem:stieltjes-carleman}
For $0<|x|<r$ and every integer $N\ge0$,
\begin{equation}
 \label{eq:stieltjes-carleman}
 |x|^{-N-1}|w(x)|
 \le C\int_{B_{3r}}|y|^{-N-1}L_r(x,y)\,|\sigma|(dy).
\end{equation}
\end{lemma}

\begin{proof}
Condition \eqref{eq:sigma-kato-newtonian} excludes atoms, so
\eqref{eq:taylor-subtracted-representation} and
\eqref{eq:sawyer-kernel-scaled} give
\[
 |x|^{-N-1}|w(x)|
 \le\int_{B_{3r}}\frac{|R_N(x,y)|}{|x|^{N+1}}\,|\sigma|(dy)
 \le C\int_{B_{3r}}|y|^{-N-1}L_r(x,y)\,|\sigma|(dy).
\]
The last integral is finite by $|\sigma|(B_\delta)=0$ and
\eqref{eq:sigma-kato-newtonian}.
\end{proof}

\subsection{Local propagation of a zero set}

Fix $r_0>0$ with $B_{8r_0}\Subset\Omega$ and take
$K_0=\overline B_{4r_0}$.  For $0<r<r_0$,
\[
 L_r(x,y)\le L_{r_0}(x,y)
 \le c_{r_0}\left(1+\log^+\frac1{|x-y|}\right).
\]
Since $B_r\subset B(y,4r)$ whenever $y\in B_{3r}$,
\begin{equation}
 \label{eq:local-kato-app}
 \sup_{y\in B_{3r}}\int_{B_r}L_r(x,y)\,\mu(dx)
 \le c_{r_0}\kappa_{K_0}(4r)\longrightarrow0
 \qquad(r\downarrow0).
\end{equation}
Choose any $r_*\in(0,r_0)$ such that
\begin{equation}
 \label{eq:kato-absorb-app}
 C|a|c_{r_0}\kappa_{K_0}(4r_*)\le\frac14.
\end{equation}
Let $u\in H^1(B_{8r_0})\cap C(B_{8r_0})$ be real-valued and satisfy
\begin{equation}
 \label{eq:measure-schrodinger}
 -\Delta u=a u\mu
\end{equation}
in distributions on $B_{8r_0}$.

\begin{lemma}
\label{lem:measure-carleman}
If $u$ vanishes on a neighborhood of the origin, then
\begin{equation}
 \label{eq:local-zero-expansion}
 u=0\qquad\text{on }B_{r_*/2}.
\end{equation}
\end{lemma}

\begin{proof}
Take $\chi\in C_c^\infty(B_{3r_*})$ with $0\le\chi\le1$ and
$\chi=1$ on $B_{2r_*}$.  Let $w,\sigma,g$ be as in
Lemma~\ref{lem:cutoff-laplacian-measure} for this cutoff.  Then
\[
 \supp\sigma\subset B_{3r_*},
 \qquad
 g\in L^2(\R^2),\qquad
 \supp g\subset B_{3r_*}\setminus B_{2r_*},
\]
and $\sigma$ satisfies \eqref{eq:sigma-kato-newtonian}.
Since $w$ vanishes near the origin, choose $0<\delta<r_*$ such that
$w=0$ on $B_\delta$.  Lemma~\ref{lem:stieltjes-carleman} then gives,
for $0<|x|<r_*$,
\begin{align}
 |x|^{-N-1}|u(x)|
 &\le C|a|\int_{B_{3r_*}}
     |y|^{-N-1}L_{r_*}(x,y)|u(y)|\,\mu(dy)\notag\\
 &\quad+C\int_{B_{3r_*}\setminus B_{2r_*}}
     |y|^{-N-1}L_{r_*}(x,y)|g(y)|\,dy.
 \label{eq:pointwise-carleman-measure}
\end{align}
For each integer $N\ge0$, define
\begin{equation}
 \label{eq:AN-def-app}
 A_N=\int_{B_{r_*}}|x|^{-N-1}|u(x)|\,\mu(dx)<\infty.
\end{equation}
Finiteness follows because $u$ vanishes near the origin.
By Tonelli and \eqref{eq:local-kato-app},
\begin{align}
 &\int_{B_{r_*}}\int_{B_{r_*}}
 |y|^{-N-1}L_{r_*}(x,y)|u(y)|\,\mu(dy)\mu(dx)\notag\\
 &\quad=
 \int_{B_{r_*}}|y|^{-N-1}|u(y)|
 \left[\int_{B_{r_*}}L_{r_*}(x,y)\,\mu(dx)\right]\mu(dy)\notag\\
 &\quad\le c_{r_0}\kappa_{K_0}(4r_*)A_N.
 \label{eq:kato-weighted-absorption}
\end{align}
For the two annular terms, put
\begin{align*}
 J_{N,1}
 &=C|a|\int_{B_{3r_*}\setminus B_{r_*}}
 |y|^{-N-1}|u(y)|
 \left[\int_{B_{r_*}}L_{r_*}(x,y)\,\mu(dx)\right]\mu(dy),\\
 J_{N,2}
 &=C\int_{B_{3r_*}\setminus B_{2r_*}}
 |y|^{-N-1}|g(y)|
 \left[\int_{B_{r_*}}L_{r_*}(x,y)\,\mu(dx)\right]dy.
\end{align*}
Integrating \eqref{eq:pointwise-carleman-measure} against $\mu$ therefore yields
\begin{equation}
 \label{eq:AN-before-remainders}
 A_N\le C|a|c_{r_0}\kappa_{K_0}(4r_*)A_N+J_{N,1}+J_{N,2}.
\end{equation}
For the first annulus, \eqref{eq:local-kato-app} bounds the inner integral.
On the second, $|x-y|\ge r_*$ gives $L_{r_*}(x,y)\le1+\log16$.
Define the finite constants
\[
 C_1=C|a|c_{r_0}\kappa_{K_0}(4r_*)
        \int_{B_{3r_*}\setminus B_{r_*}}|u|\,d\mu,
 \qquad
 C_2=C(1+\log16)\mu(B_{r_*})\|g\|_1.
\]
They are independent of $N$, and
\begin{equation}
 \label{eq:annular-remainders-app}
 J_{N,1}\le C_1r_*^{-N-1},
 \qquad
 J_{N,2}\le C_2(2r_*)^{-N-1}.
\end{equation}
Together with \eqref{eq:kato-absorb-app}, this gives
\begin{equation}
 \label{eq:AN-final-app}
 A_N\le\tfrac14 A_N+(C_1+C_2)r_*^{-N-1}
 \quad\Longrightarrow\quad
 A_N\le\tfrac43(C_1+C_2)r_*^{-N-1}.
\end{equation}
For every $0<\theta<1$,
\[
 \int_{B_{\theta r_*}}|u|\,d\mu
 \le(\theta r_*)^{N+1}A_N
 \le\tfrac43(C_1+C_2)\theta^{N+1}
 \longrightarrow0
 \qquad(N\to\infty).
\]
Consequently, on $B_{\theta r_*}$,
\[
 u=0\quad\mu\text{-almost everywhere},
 \qquad
 -\Delta u=a u\mu=0\quad\text{in distributions}.
\]
The function $u$ is therefore harmonic and real analytic there.
Its vanishing near the origin gives
\[
 u=0\qquad\text{on }B_{\theta r_*}.
\]
Taking $\theta=1/2$ proves
\eqref{eq:local-zero-expansion}.
\end{proof}

\subsection{Global continuation}

Assume $\Omega$ is connected.  Let
$u\in H^1_{\mathrm{loc}}(\Omega)\cap C(\Omega)$ be real-valued and satisfy
\begin{equation}
 \label{eq:global-measure-eigen}
 -\Delta u=a u\mu
\end{equation}
in distributions on $\Omega$, and let $U\subset\Omega$ be nonempty and open.

\begin{theorem}
\label{thm:measure-ucp}
We have
\begin{equation}
 \label{eq:global-ucp-conclusion}
 u=0\text{ on }U
 \quad\Longrightarrow\quad
 u\equiv0\text{ in }\Omega.
\end{equation}
\end{theorem}

\begin{proof}
Suppose $u=0$ on $U$.  Fix $p\in U$ and $q\in\Omega$, and choose a
polygonal path $\Gamma:[0,1]\to\Omega$ from $p$ to $q$.
Its image $K=\Gamma([0,1])$ is compact.  Choose $r_0>0$ so that
the closed $8r_0$-neighborhood of $K$ lies in $\Omega$, and put
\[
 K_{4r_0}=\{z\in\R^2:\operatorname{dist}(z,K)\le4r_0\}.
\]
By \eqref{eq:measure-kato-global}, there is $0<r_K<r_0$ such that
\begin{equation}
 \label{eq:uniform-kato-absorption}
 C|a|c_{r_0}\sup_{x\in K}
     \kappa_{\overline{B(x,4r_0)}}(4r_K)
 \le C|a|c_{r_0}\kappa_{K_{4r_0}}(4r_K)
 \le\frac14.
\end{equation}
Lemma~\ref{lem:measure-carleman}, after translation, therefore applies with
$r_*=r_K$ at every $x\in K$.
Choose successive points $x_0=p,\ldots,x_m=q$ on $\Gamma$ with
$|x_{j+1}-x_j|<r_K/4$.  Then
\[
 B(x_{j+1},r_K/4)\subset B(x_j,r_K/2),
 \qquad 0\le j<m,
\]
and the local result gives the induction
\[
 u=0\text{ near }x_j
 \ \Longrightarrow\
 u=0\text{ on }B(x_j,r_K/2)
 \ \Longrightarrow\
 u=0\text{ near }x_{j+1}.
\]
Starting at $p\in U$ yields $u(q)=0$.  Since $q\in\Omega$ was arbitrary,
\eqref{eq:global-ucp-conclusion} follows.
\end{proof}

\Needspace{11\baselineskip}
On the probability-one event of Section~\ref{sec:prelim}, let $\lambda>0$
and let $u\in H_0^1(D)\cap C(D)$ be real-valued, with
\begin{equation}
 \label{eq:lbm-eigenfunction-measure-equation}
 \E(u,v)=\lambda\int_Duv\,dM_h,
 \qquad v\in C_c^\infty(D).
\end{equation}

\begin{corollary}
\label{cor:lbm-eigenfunction-ucp}
We have
\begin{equation}
 \label{eq:lbm-eigenfunction-ucp}
 u=0\text{ on a nonempty open subset of }D
 \quad\Longrightarrow\quad
 u\equiv0\text{ in }D.
\end{equation}
\end{corollary}

\begin{proof}
By the normalization of $\E$,
\[
 \langle-\Delta u,v\rangle
 =2\pi\E(u,v)
 =2\pi\lambda\int_Duv\,dM_h,
 \qquad v\in C_c^\infty(D).
\]
Thus $-\Delta u=2\pi\lambda uM_h$ in distributions.
The domain $D$ is connected, $u\in H^1_{\mathrm{loc}}(D)\cap C(D)$, and
the locally finite measure $M_h$ satisfies \eqref{eq:measure-kato-global}
by \eqref{eq:kato-gmc}.  Theorem~\ref{thm:measure-ucp} applies with
$a=2\pi\lambda$.
\end{proof}

\section*{Use of AI-assisted tools}
ChatGPT (OpenAI) was used to assist with checking journal submission requirements,
verifying bibliographic details, and preparing editorial revisions. The author
is responsible for reviewing the changes and for the final content of the manuscript.

\section*{Statements and Declarations}
\noindent\textbf{Competing interests.} The author declares no competing interests.\par
\noindent\textbf{Data availability.} No datasets were generated or analyzed during
the current study.

\end{document}